\documentclass{article}
\usepackage{amssymb}
\usepackage{latexsym}
\newtheorem{thm}{Theorem}
\newtheorem{tom}{Theorem}
\newtheorem{lem}{Lemma}

\newtheorem{dfn}{Definition}
\begin{document}
\title{Calibrated Fibrations on Complete Manifolds via Torus Action}
\author{Edward Goldstein}
\maketitle

\renewcommand{\abstractname}{Abstract}
\begin{abstract}
In this paper we will investigate torus actions on complete manifolds with
calibrations. For Calabi-Yau manifolds $M^{2n}$ with a
Hamiltonian structure-preserving k-torus action we show that 
any smooth symplectic reduction has a natural holomorphic volume form. Moreover
Special Lagrangian (SLag) submanifolds of the reduction lift to SLag submanifolds of $M$, invariant under the torus action.
If $k=n-1$ and $H^1(M,\mathbb{R})= 0$ then we prove that $M$ is a
fibration with generic fiber being a SLag submanifold.
As an application we will see that crepant  resolutions of singularities of a
finite Abelian subgroup of
$SU(n)$ acting on $\mathbb{C}^n$ have SLag fibrations.
We study SLag submanifolds on the total space $K(N)$ of the canonical bundle of a Kahler-Einstein manifold $N$ with positive scalar curvature and give a conjecture about fibration of $K(N)$ by SLag sub-varieties, which we prove if $N$ is toric. We will also get somewhat weaker results for coassociative submanifolds of a $G_2$-manifold $M^7$, which admits a 3-torus, a 2-torus or an
$SO(3)$ action.
\end{abstract}

\section{Introduction}
In this paper we will use structure preserving torus actions on non-compact
manifolds with calibrations to construct calibrated submanifolds (both for
Calabi-Yau manifolds and for 7-manifolds with a $G_2$ structure). We will assume that no element of the torus acts trivially. Throughout the paper we will use the notion of a calibrated fibration: 
\begin{dfn}
 Let $(M,\varphi)$ be a Riemannian manifold with a calibrating form $\varphi$. Then we say that $M$ has a calibrated fibration on it if there is a surjective map $\alpha: M \mapsto V$ onto a topological space $V$ and a subset $S \subset M$ s.t.
 
i) For any point $m \in M-S$ the level set $L_m$ of $\alpha$ through $m$ is a smooth submanifold, calibrated by $\varphi$. 

ii) The set $S$ is locally contained in a finite union of submanifolds of codimension $\geq 4$ in $M$.
\end{dfn}
In Section 2 we consider a Kahler manifold $(M^{2n},\omega)$ with a
non-vanishing holomorphic $(n,0)$ form $\varphi$. We can define, as in
\cite{Gold},  Special
Lagrangian (SLag) submanifolds $L$ by the conditions:
 \[\omega |_L=0 ~ , ~ Im\varphi |_L=0.\] 
If $g$ is the Riemannian metric corresponding to $\omega$, then we can
conformally scale $g$ to a metric $g'$ on $M$ so that the form $\varphi$ will
have length $\sqrt{2}^n$ with respect to $g'$. Then SLag
submanifolds will be calibrated by $Re\varphi$ with respect to $g'$. In
particular, they will be minimal submanifolds of $(M,g')$ and Lagrangian
submanifolds of $\omega$.

If $M$ is compact and simply connected, then for any Kahler form $\omega$ on $M$ Yau's celebrated resolution of the Calabi conjecture gives a (unique) Ricci-flat Kahler form $\omega'$ in the same cohomology class as $\omega$ (see \cite{Yau}). Also the SYZ conjecture (see \cite{SYZ}) states that
$(M,Re \varphi)$ has a calibrated fibration
with generic fiber being a SLag torus with respect to a Ricci-flat Kahler metric. We can ask an analogous question for any Kahler metric on $M$ and
we showed in \cite{Gold} that this holds for a choice of Kahler metric on
a Borcea-Voisin threefold. In this paper we will be interested in non-compact
Calabi-Yau manifolds with a structure-preserving torus action. The main results of Section 2 are as follows:
\begin{tom}
Suppose we have a Hamiltonian structure-preserving
$T^k$-action on $M^{2n}$. Then any smooth symplectic reduction $M_{red}$ has a
natural holomorphic volume form $\varphi_{red}$. Moreover SLag
submanifolds of $M_{red},\omega_{red},\varphi_{red}$ lift to SLag submanifolds of $M$, invariant under $T^k$-action. Vice versa, let $L$ be a connected, $T$-invariant SLag submanifold of $M$ s.t. $T$ acts freely on $L$. Then $L$ lies on a level set $a$ of the moment map $\mu$ and one can find a SLag submanifold $L'$ in the smooth part of the symplectic reduction through level set $a$ s.t. $L'$ lifts to $L$.
\end{tom}
\begin{tom}
 Suppose that $k=n-1$ and $H^1(M,\mathbb{R})=0$. Then $M$ has a calibrated fibration $\alpha$ over an open subset of $\mathbb{R}^n$ with the set $S$ of singular points being the non-regular points of the torus action (i.e points there the differential of the action is not injective). Moreover for a generic point $p$ (outside of a countable union of $(n-2)$-dimensional planes in $\mathbb{R}^n$), the fiber $\alpha^{-1}(p)$ is a smooth SLag submanifold. Connected components of each smooth fiber are diffeomorphic to an
$(n-1)$-torus times $\mathbb{R}$. Singular fibers have singularities of codimension at least 2, and near singular points they are diffeomorphic to a product of a cone with a Euclidean ball.
\end{tom}
If we make certain assumptions on the set of non-regular points of $T$-action, then we can replace the countable union by a finite union in Theorem 2 (see Theorem 3 in Section 2).

In Section 3 we consider a 7-manifold $M$ with a $G_2$-form $\varphi$
(see \cite{BS}). 
Let $H^1(M,\mathbb{R})=0$. If $M$ has a 3-torus action, then $M$ is covered by a family of non-intersecting coassociative submanifolds. Suppose
$M$ admits a 2-torus action. Then we will define certain $G_2$-reductions
$M_{red}$,
which are symplectic 4-manifolds with a compatible almost complex structure and trivial canonical bundle. We
will see that 2-dimensional 
complex sub-varieties of $M_{red}$ lift to $T$-invariant, coassociative submanifolds of $M$. Finally suppose $M$ admits an $SO(3)$-action, which is not regular in at least one point. If $SO(3)$ acts freely on the set $M'$ of regular
points of the $SO(3)$-action, then $M'$ is covered by a family of non-intersecting coassociative submanifolds, diffeomorphic to $SO(3) \times \mathbb{R}$.

In Section 4 we will consider some applications of the results of the two previous
sections. For the Calabi-Yau case, we will show that crepant resolutions of singularities of a finite Abelian subgroup of $SU(n)$ acting on $\mathbb{C}^n$ have SLag fibrations. Also for any Kahler-Einstein manifold $N$ with positive scalar curvature, the total space $K(N)$ of it's canonical bundle is a Calabi-Yau manifold (see \cite{Sol}).
We  investigate SLag submanifolds on $K(N)$. For each orientable minimal Lagrangian submanifold of $N$ we associate a 1-parameter family $(L_{\lambda}|\lambda \in \mathbb{R})$ of SLag submanifolds of $K(N)$. Also $L_0$ is invariant under scaling of $K(N)$ by a real number. For any compact Kahler manifold $N^{2n}$ with an effective n-torus action we prove that one of the regular orbits of the action is a minimal Lagrangian submanifold of $N$. If $N$ is Kahler-Einstein with nonzero scalar curvature $t$ and toric, then we prove that precisely 1 such orbit $L$ is a minimal Lagrangian submanifold of $N$. For $t>0$ we use Theorem 2 to construct a SLag fibration on $K(N)$ and we prove that all fibers are asymptotic at infinity to the fiber $L_0$.
We conjecture that any K-E manifold $N$ with positive scalar curvature has a minimal Lagrangian submanifold $L$. Moreover $K(N)$ fibers with generic fiber being a SLag submanifold of $K(N)$ and all fibers are asymptotic to $L_0$ at infinity.

In the $G_2$ case, Bryant and Salamon have constructed in \cite{BS} some examples of complete metrics with holonomy
$G_2$. Some metrics are on the total space of the spin bundle over a
3-dimensional space form. Others are on a total space of a bundle
$\Lambda^2_-$ of anti-self-dual 2-forms over a self-dual Einstein
4-manifold. Many examples admit $T^2$ and $SO(3)$-actions, and we show that in one example the $G_2$ manifold $M$ can be covered by non-intersecting coassociative $SO(3)$-invariant submanifolds.
  
{\bf Acknowledgments} : This paper is a part of author's work towards his
Ph.D. at the Massachusetts Institute of Technology. The author wants to
express his gratitude to his advisor, Tom Mrowka, for initiating him into
this subject and for continuing support.

After writing this paper the author learned that Mark Gross has independently obtained results, which are similar to some results of Sections 2  and 4.3.
\section {Torus action on Calabi-Yau manifolds}
Let $(M^{2n},\omega,\varphi)$ be a Calabi-Yau manifold with a structure-preserving
Hamiltonian $T^k$-action. For any element $v$ of the Lie algebra
$\mathcal{G}$ of $T^k$ we
associate the infinitesimal flow vector field $X_v$ on $M$, induced by
the differential of the action. Then the $X_v$ commute and their flows
preserve $\omega$ and $\varphi$. Let $v_1, \ldots, v_l$ be elements of
${\cal G}$ and $X_1, \ldots, X_l$ be the corresponding vector fields on
$M$. We claim that the $(n-l,0)$- form $\varphi'= i_{X_1} \ldots i_{X_l}
\varphi$ (obtained by contraction of $\varphi$ by the vector fields $X_1,
\ldots, X_l$) is a closed $(n-l,0)-$form. We will prove it by induction on
$l$. If $l=1$ then the $X_1$-flow preserves $\varphi$ and so \[ 0 =
{\cal L}_{X_1}\varphi= d(i_{X_1}\varphi)=d\varphi'\]
Now we use induction on $l$. The $X_1$ flow preserves $X_2, \ldots, X_l$
and it preserves $\varphi$, so it preserves $\varphi''=i_{X_2} \ldots
i_{X_l}\varphi$. So \[0 = {\cal L}_{X_1}\varphi''=
d(i_{X_1}\varphi'') + i_{X_1}(d\varphi'') \] and we are done by induction.

A moment map $\mu$ for the $T^k$-action is a map $\mu: M \mapsto
{\cal G}^{\ast}$ (the dual Lie Algebra of $T^k$), which satisfies 
\[d(\mu(v))= i_{X_v}\omega \] for any $v \in \mathcal{G}$. The moment map is $T^k$-invariant. 
\begin{thm}
Any smooth symplectic reduction $M_{red}$ has a natural holomorphic volume form
$\varphi_{red}$. Moreover SLag submanifolds of $M_{red}$
lift to SLag submanifolds of $M$, invariant under
$T^k-$action. Vice versa let $L$ be a connected, $T$-invariant SLag submanifold of $M$. Then $L$ lies on a level set $a$ of the moment map $\mu$ and moreover there is a SLag submanifold $L'$ on a smooth part of symplectic reduction through level set $a$ s.t. $L$ is the lift of $L'$.
\end{thm}
{\bf Proof : } Consider a level set $\Sigma_a$ of the moment map
$\mu$, on which $T^k$ acts freely. Since
$\mu$ is $T^k-$invariant, the vector fields $X_v$ are tangent to
$\Sigma_a$. Consider the bundle $V= ~span(X_v)$ over $\Sigma_a$. Since $d(\mu(v))= i_{X_v}\omega$, the
tangent bundle to
$\Sigma_a$ is the $\omega-$orthogonal complement of $V$. Also $V
\bigcap JV = 0$ (here $J$ is the complex structure on $M$).  
Let $W = (V \oplus JV)^{\perp}$ ( here $\perp$ is with respect to the
metric). Then $W$ is a complex vector bundle over $\Sigma_a$, the tangent
bundle to $\Sigma_a$ is $W \oplus V$ and the quotient of $W$ by the $T$-action can be viewed as a tangent
bundle to the symplectic reduction $M_{red}= \Sigma_a/ T^k$. 

Let $v_1, \ldots, v_k$ be a basis for the Lie algebra of $T^k$ and $X_1,
\ldots, X_k$ corresponding vector fields on $M$. Let $\varphi'= i_{X_1}
\ldots i_{X_k} \varphi$. Then as we have seen, $\varphi'$ is a holomorphic
$(n-k,0)$ form on $M$. Also $\varphi'|_W$ is a holomorphic volume form on $W$.
Since $\varphi'$ is $T^k$-invariant, it is clear that there is a unique
$(n-k,0)$ form $\varphi_{red}$ on $M_{red}$
s.t. $\pi^{\ast}(\varphi_{red})=\varphi'$ on $\Sigma_a$ (here
$\pi: \Sigma_a \mapsto M_{red}$ is the quotient map). 
Now $\varphi_{red}$ is closed (since $\varphi'$ is), and hence it is a holomorphic volume form on $M_{red}$.

Let $L'$ be a SLag submanifold of $M_{red}$ and $L=
\pi^{-1}(L')$. It is clear that $L$ is a SLag submanifold of
$M$, invariant under $T^k$-action. Vice versa, let $L$ be a connected
SLag submanifold of $M$, invariant under $T^k$-action and $T$ acts freely on $L$. Since
$L$ is invariant under the torus action, the vector fields $X_v$ are
tangent to $L$. Let $u$ be some vector, tangent to $L$. Since $L$ is
Lagrangian, we have \[ 0 = \omega(X_v,u)= d(\mu(v))(u) \]
So the differential of the moment map is zero on $L$, hence $L$ lies on some
level set $\Sigma_a$ of the moment map. Since $T$-action is free on $L$, it is freely on a neighbourhood $U$ of $L$ in $\Sigma_a$. Moreover $U$ is a smooth submanifold of $M$. So $U_{red}=U/T$ will be an open set in the smooth part of the symplectic reduction through $a$.
It is clear that the quotient $L'= L / T^k$ is a SLag submanifold on $U_{red}$.   Q.E.D.

In case $k=n-1$ $M$ has a calibrated fibration by SLag submanifolds, invariant under $T^k$-action:
\begin{thm} Let $k=n-1$ and $H^1(M,\mathbb{R})=0$. Then 

i) $M$ has a calibrated fibration $\alpha$ over an open subset of $\mathbb{R}^n$ with the set $S$ of singular points being the non-regular points of the $T$-action. 

ii) For a  generic point $p$ (outside of a countable union of $(n-2)$-planes in $\mathbb{R}^n$), the fiber $\alpha^{-1}(p)$ is a smooth SLag submanifold of $M$.

iii) Connected components of each smooth fiber are diffeomorphic to a product of an $(n-1)$-torus with $\mathbb{R}$.

iv) Singular fibers have singularities of codimension at least 2, and near a singular point they are diffeomorphic to a product of a cone with a Euclidean ball. 
\end{thm}
{\bf Proof:} Define the form $\varphi'$ as in the proof of Theorem 1. Then
$\varphi'$ is a holomorphic $(1,0)$-form, invariant under the torus action. Since $H^1(M,\mathbb{R})=0$, there is a holomorphic function $f=\eta+ i\xi$ s.t. $df=\varphi'$. It is
clear that $f$ is also invariant under $T^k$-action. Let $L$ be a connected SLag submanifold, invariant under the torus action. As we have seen,
$L$ must lie on the level set of the moment map $\mu$. Also, since $L$ is
Special, one easily deduces that $Im\varphi'|_{L}=0$, so $L$ must lie on
a level set of $\xi=Imf$, i.e. $L$ lies on a level set of n-functions \[
\mu= a ~ , ~ \xi=c \]
The moment map $\mu$ goes to ${\cal G}^{\ast}$, which we identify with $\mathbb{R}^{n-1}$ by choosing a basis of ${\cal G}$. We define $\alpha=(\mu,\xi): M \mapsto \mathbb{R}^n$.

Let $S$ be the set of non-regular points of the torus action.
We claim that a level set $L_m$ of $(\mu,\xi)$, that passes through a regular point $m \in M-S$ is smooth n-dimensional SLag submanifold of $M$ near $m$.
Indeed let $\Sigma_a$ be the level set of the moment map passing through $m$
and $V$ and $W$ be vector bundles on $\Sigma_a$ near $m$ as in the proof of Theorem 1. Let $v_1, \ldots, v_k$ be
a basis for the Lie algebra of $T^k$. Then $d\mu(v_1), \ldots, d\mu(v_k)$
is basis of  $(JV)^{\ast}$. Also those 1-forms vanish on $W$. Now
$d\xi=Im\varphi'$ restricted to $W$ is non-zero. So the forms
$d\mu(v_1), \ldots, d\mu(v_k), d\xi$ are linearly independent at $m$, and so
the level set $L_m$ is a smooth submanifold of $M$ near $m$.

Next we prove that $L_m$ is SLag near $m$. Since $\xi$ and $\mu$ are
$T^k$-invariant, then so is $L_m$. So the $X_v$ are in the tangent space to $L_m$ at $m$.
Since $L_m$ is on the level set of $\mu$, the tangent space to $L_m$ at $m$ is
$\omega$-orthogonal to $X_v's$, so it must be Lagrangian. Also
$Im\varphi'|_{L_m}=0$ implies that $L_m$ is Special at $m$.

Now we prove {\bf iii)}: Let $L$ be a level set of $(\mu,\xi)$ s.t. all points on $L$ are regular points for $T$-action on $M$.
Let $L'$ be a connected component of $L$. We have $Re\varphi'=d\eta$. 
One easily sees that $Re\varphi'|_{L} \neq 0$ for all points of $L$. One also easily shows that $\nabla \eta$ along $L$ is tangent to $L$, hence it coincides with $\nabla (\eta|_L)$.

 Consider a $T^{n-1}$-orbit $T$ on $L'$. $T$ lives on some level set of $\eta$ on $L'$, hence it must coincide with a connected component of this level set. The normalized gradient flow of $T$
by $\frac{\nabla \eta}{|\nabla \eta|^2}$ on $L'$ gives a diffeomorphism between $L'$ and $T^{(n-1)}$ times $\mathbb{R}$. Indeed this flow is defined on an interval $(a_-,a_+)$, there $a_-,a_+$ are independent on the choice of a point in $T$. One easily deduces that the orbit of $T$ under this flow is isolated in $L$, hence it is equal to $L'$.  

Next we prove {\bf ii)} and show that $S$ is locally contained in a finite union of submanifolds of codimension 4 in $M$. 
To prove this we need to understand the picture near a point in $S$.
Let $m \in S$. The differential of the action is not injective at $m$ and $m$ has a stabilizer $T'$ of positive dimension $l$ and an orbit $O$. To prove {\bf ii)} we need to see what is the image of $(\mu,\xi)$ on $S$ near $m$. 

The symplectic form $\omega$ restricts trivially to $O$, hence we have $\omega=d\gamma$ for some 1-form $\gamma$ in a neighbourhood of $O$. Since $\omega$ is $T^{n-1}$-invariant, we can make $\gamma$ invariant as well (by integrating over $T^{n-1}$). For any $v \in {\cal G}$ we have $0 = {\cal L}_{X_v}\gamma= d(\gamma(X_v))+i_{X_v}\omega$. So the map $\mu'(v)=-\gamma(X_v)$ is a moment map near $O$ and $\mu-\mu'$ is a constant. Obviously $Im\varphi|_O=0$, so we can write $Im\varphi=d\beta$ for a $T$-invariant $(n-1)$-form $\beta$ in a neighbourhood of $O$. Let $\xi'=\beta(X_1,\ldots,X_{n-1})$. Arguing as in the proof of Theorem 1 we get that $\varphi'=d\xi'$, so $\xi-\xi'$ is a constant (here $\xi=Imf$ as before). We will prove that the image of $S$ near $O$ by $(\mu',\xi')$ is contained in a finite union of $(n-2)$-planes in $\mathbb{R}^n$. This will prove that one can find a neighbourhood $U$ of $m$ and a finite union $H$ of $(n-2)$-planes in $\mathbb{R}^n$ s.t $(\mu,\xi)(S \bigcap U) \subset H$. Since $M$ is paracompact, one can find a countable union of $(n-2)$-panes $H'$ in $\mathbb{R}^n$ s.t. $(\mu,\xi)(S) \subset H'$ and ii) follows. 

Obviously $\xi'=\beta(X_1,\ldots,X_{n-1})=0$ on $S$. Next we prove that there is a collection $v_1,\ldots,v_l$ of linearly independent elements of ${\cal G}$ s.t. at any point $p' \in S$ the flow field corresponding to $v_i-v_j$ vanishes for some $i$ and $j$. This will imply that the image of $\mu'$ on $S$ will lie on a finite collection of hyper-planes in the dual Lie Algebra ${\cal G}^{\ast}$. Hence the image of $(\mu',\xi')$ on $S$ near $O$ will be contained in a finite collection of $(n-2)$-planes. 
 
We have the tangent bundle $TO$ and the normal bundle $N(O)$, which splits as a direct sum $N(O)=J(TO) \oplus W$, $W=(TO+J(TO))^{\perp}$. $W$ is a complex vector bundle of dimension $l+1$. $T'$-action preserves $O$ and it's differential preserves $TO$ and $J(TO)$, hence $T'$ acts faithfully on $W$ ( because no element of $T'$ acts trivially on $M$). We get an injective homomorphism $\rho$ from $T'$ to $SU(W)$, whose image is in some maximal torus of $SU(W)$. By dimension count this image is the maximal torus of $SU(W)$. We identify a small neighbourhood $V$ of $O$ in $M$ with a small ball in $N(O)$ by the exponential map. Then the action of $T'$ under this identification is trivial on $J(TO)$ and equal to the action by representation $\rho$ on $W$. Let $\pi: V \simeq N(O) \mapsto O$ be the projection. Take any element $v$ of the Lie Algebra of $T$. Then it is clear that $\pi_{\ast}(X_v)=0$ iff $v$ is in the Lie Algebra ${\cal G}'$ of $T'$. So a point $p \in V$ is in $S$ iff there is an element $0 \neq v \in {\cal G}'$ s.t. $X_v(p)=0$.

We identify the fiber $W(m)$ at $m$ with $\mathbb{C}^{l+1}$. We have
a torus $T' \simeq T^l$ acting as a maximal torus of $SU(l+1)$ on $\mathbb{C}^{l+1}$ by 
\begin{equation}
(e^{i\theta_1},\ldots,e^{i\theta_l})(z_1,\ldots,z_{l+1})=(e^{i\theta_1}z_1,
\ldots,e^{\theta_l}z_l,e^{-\Sigma \theta_i}z_{l+1}) 
\end{equation}
It is clear that the the non-regular points in $\pi^{-1}(m)$ are subspaces $H_{i,j}=(z_i=z_j=0) \oplus J(TO)$ and the vector field $\partial_{\theta_i}-\partial_{\theta_j}$ is in the kernel of the differential of the action at these points. Let $U'$ be some submanifold of $T$, which passes through $0 \in T$ and is transversal to $T'$ at $0$.
Then the image of $m$ under $U'$-action is a neighbourhood of $m$ in $O$. Let $U_{i,j}$ be the image of $H_{i,j}$ under $U'$-action. The flow vector field of $\partial_{\theta_i}-\partial_{\theta_j}$ vanishes along $U_{i,j}$. We claim that near $m$ $S$ is contained in the union of $U_{i,j}$.
Indeed let $p' \in S$ and let $m'=\pi(p')$. There is an element of $U'$ which sends $m$ to $m''$ and hence it sends $p$ to $p'$ for some $p \in \pi^{-1}(m)$.  One easily deduces that $p \in \bigcup H_{i,j}$, hence $p' \in \bigcup U_{i,j}$. Thus we can take $v_i=\partial_{\theta_i}$ and we are done. Moreover $U_{i,j}$ is obviously a submanifold of $M$ of codimension 4, thus we also proved that $S$ is locally contained in a finite union of submanifolds of codimension 4 in $M$.

To complete the proof of {\bf i) and iv)} we still need to investigate the structure of the singular fiber $L_m$ through $m$ and to prove that the image of $(\mu,\xi)$ is open. Let $e_1,\ldots,e_l$ be the basis of Lie algebra of $T'$. We extend it by $e_{l+1},\ldots,e_{n-1}$ to be the basis of ${\cal G}$. Let $\mu''= (\mu(e_{l+1}),\ldots,\mu(e_{n-1}))$. Then the differential of $\mu''$ is surjective along $O$ and the level set $\Sigma$ of $\mu''$ through $m$ is a smooth submanifold of $M$ (containing $O$). Also obviously $L_m \subset \Sigma$. We can investigate $L_m$ by means of a local symplectic reduction.
Let $Q=span(e_{l+1},\ldots,e_{n-1})$. Take a small ball $U$ containing the origin in $Q$. $U$ can be identified with a submanifold (still called $U$) of $T^{n-1}$ via the exponential map. Also consider the induced metric on $\Sigma$ and let $Z$ be the image of a small ball in the normal bundle to $O$ in $\Sigma$ at $m$ by the exponential map. So $Z$ is $T'$-invariant, contains $m$ and is transversal to $O$. We will define an equivalence relation on a small neighbourhood $V'$ of $m$ in $\Sigma$ by making the equivalence classes to be the orbits of $U$-action through points of $Z$. The quotient $M'$ can be thought of a local symplectic reduction of $M$ by the action of $U$. So $M'$ is a Kahler manifold. By Theorem 1 we have a natural trivialization $\varphi''$ of the canonical bundle of $M'$. We have a structure-preserving $T'$-action on $M'$. Let $\mu^{\ast}$ be the restriction of $\mu$ to the dual Lie Algebra of $T'$. Then $\mu^{\ast}$ is a moment map for $T'$-action on $M'$. Also $\xi$ descends to $M'$ and the level sets of $(\mu^{\ast},\xi)$ are SLag submanifolds of $M'$, which lift to level sets of $(\mu,\xi)$ on $M$. 

We will investigate the level sets of $(\mu^{\ast},\xi)$ on $M'$. Let $\tau :V' \mapsto M'$ be the quotient map. Then $\tau: Z \mapsto M'$ is a diffeomorphism.
Let $L'_m$ be a level set of $(\mu^{\ast},\xi)$ through $\tau(m)$. We will prove that $L'_m$ is diffeomorphic to an $(l+1)$-dimensional cone and moreover all points on $L'_m - \tau(m)$ are regular points for the $T'$-action. We claim that {\bf  i) and iv)} follow from this. Indeed $L_m$ is an orbit of the $U$-action on $\tau^{-1}(L'_m) \bigcap Z$. Thus $L_m$ will be locally a product of a cone with a Euclidean ball. Also it's singular set is of codimension $l+1 \geq 2$.
Moreover $L_m$ will contain regular points for the $T^{n-1}$-action. The differential of $(\mu,\xi)$ is surjective at those points, and thus we will deduce that the image of $(\mu,\xi)$ is open.

So we have a Kahler manifold $M'$ with a trivialization $\varphi''$ of the canonical bundle and a structure-preserving $T' \simeq T^l$-action, which preserves $m'=\tau(m) \in M'$ and induces an action of a maximal torus of $SU(l+1)$ on the tangent space at $m'$ as in equation (1). We would like to understand the level set $\mu^{\ast}(m')$ of the moment map $\mu^{\ast}$ (which contains $L'_m$). To do this we introduce Equivariant Darboux coordinates in a following way : We identify a small neighbourhood of $M'$ with a ball $Y$ on a tangent bundle $T_{m'}M'$ via the exponential map. The action of $T'$ will be linear on $Y$, and it will preserve the symplectic form $\omega'$ on $Y$, induced by the exponential map. We identify $Y$ with a ball in $\mathbb{C}^{l+1}$. Then we will also have a standard symplectic form $\omega_0$ on $Y$. Moser's proof of the Darboux theorem gives a embedding $\phi$ of a possibly smaller ball $Y'$ into $Y$, s.t $\phi^{\ast}(\omega_0)=\omega'$. Now both $\omega'$ and $\omega_0$ are $T'$-invariant, so $\phi$ will be $T'$-equivariant. 

The moment map of $\omega^0$ is $\mu^0= (\mu_1,\ldots,\mu_l)$ with $\mu_i=|z_i|^2-|z_{l+1}|^2$. The non-regular points of the action are, as we saw, the planes $(z_i=z_j=0)$ and the zero set of $\mu^0$ intersects them only at the origin.
The zero set of $\mu^0$ is a cone $|z_i|=|z_j|$ in $\mathbb{C}^{l+1}$. 
It's symplectic reduction thus will be a 2-dimensional cone with a singular point at the origin. Similarly we can take a level set $P_0$ of $\mu^{\ast}$ through $m'$ to get a symplectic reduction $M''$. Let $m'' \in M''$ be the image of $m'$ under the quotient map.
We take a fixed compact neighbourhood $K$ of $m''$ in $M''$ and $K$ will be a 2-dimensional cone with boundary.

We had a holomorphic function $\eta+i\xi$ on $M$ as before, and this function descends to a holomorphic function on $M'$ and on $M''$. W.l.o.g. we assume that $\xi(m)=0$.
The zero set of $\xi$ on $M''$ lifts to the fiber $L'_m$. The gradients of $\xi$ and $\eta$ on $M$ are orthogonal to the orbits of $T^{n-1}$-action. Thus the gradient flow of, say, $\eta$ on $M$ projects to the gradient flow of $\eta$ on $M'$ and on $M''$. 
Also the gradients of $\xi$ and $\eta$ are linearly independent over $M''- m''$ and $\nabla\eta = J\nabla \xi$. Thus the gradient flow of $\eta$ is the Hamiltonian flow of $\xi$, and so it preserves $\xi$. Since $d\xi \neq 0$ on $M''- m''$ then near every point on the zero set of $\xi$ in $M''-m''$, $\xi^{-1}(0)$ is an orbit of the gradient flow of $\eta$. The orbits of this gradient flow on $\xi^{-1}(0)$ are obviously isolated. One end of such orbit might flow to the critical point $m''$, but the other must flow to the boundary of $K$. From this we deduce that there are only finitely many of these orbits in $\xi^{-1}(0)$ in $K$. We will look even for a smaller $K' \subset K$ so that one end of each orbit in $K'$ will flow to $m''$. So these will be orbits $d_1,\ldots, d_p$. We will prove that $p > 0$. We claim that from this it follows that $L'_m$  diffeomorphic to  a cone modeled on a $p$ $(l)$-tori (i.e $p$ $l$-tori will be the base of the cone).

Indeed let $P_0$ be the level set $\mu^{\ast}(m')$ of $\mu^{\ast}$ on $M'$ as before and let
$\nu :P_0 \mapsto M''$ be the quotient map. Take a point $q_i$ on $\nu^{-1}(d_i)$  in $M'$. Then the gradient flow of $\eta$ through $q_i$ projects under $\nu$ to the gradient flow on $M''$. Hence the gradient flow of $\eta$ through $q_i$ must terminate at $m'$. Let $d'_i$ be the trajectory of this flow. Then the orbit $D_i$ of $T'$-action on $d'_i$ is diffeomorphic to a cone modeled on an $l$-torus. Moreover $\nu^{-1}(d_i)=D_i$. So $L'_m$ is diffeomorphic to a cone, modeled at $p$ $l$-tori.   

Finally we prove that $p>0$.
We have $\mu^0$ and on every non-regular point of $T'$-action some of the functions $\mu_i-\mu_j$ vanish. The set of all such points on which some of $\mu_i-\mu_j$ vanish is a union of hypersurfaces in $M'$. If we take a point $s$ outside of these hypersurfaces then the level set $P_s$ of the moment map $\mu^{\ast}$ containing $s$ is smooth. Take now a sequence of such points $s_j$ converging to $m'$. Fix a compact neighbourhood $B$ of $m'$ in $M'$. We consider the positive direction $\eta$-gradient flow  lines $a_j$ through $s_j$. There are no critical points of $\eta$ on $P_{s_j}$. So those lines $a_j$ must intersect the boundary $\partial B$. We saw that the level set of $P_0$ is smooth outside of $m'$. So $P_{s_j}$ converge to $P_0$ outside of $m'$.
It is easy to show that (after passing to a subsequence) $a_j$ converge to a flow line on $P_0$ terminating at $m'$. We project it to $M''$ and get a gradient flow line on $M''$ terminating at $m''$ and we are done.   Q.E.D. 

{\bf Remark}: The proof of Theorem 2 in fact gave a construction of Special
Lagrangian submanifolds as level sets \[ \mu=a ~ , ~ \xi=c\] Thus we
effectively got an algebraic construction of Special Lagrangian
submanifolds. We will  utilize this construction for some examples in Section 4.   

The countable union of planes in Theorem 2 stems from the fact that $M$ is non-compact. If we make certain assumptions on the set of non-regular points of the $T$-action, then we can replace the countable union by a finite union:
\begin{thm}
Let $k=n-1$ as in Theorem 2. Suppose that the set of non-regular points of the $T$-action on $M$ is a finite union $S= \bigcup S_i$ of connected submanifolds s.t. each $S_i$ has a positive-dimensional stabilizer $T_i \subset T$. Then for all points $p$ outside a finite union $H= \bigcup H_i$ of $(n-2)$-planes in $\mathbb{R}^n$, the fiber $(\mu,\xi)^{-1}(p)$ is a smooth SLag submanifold of $M$.
\end{thm}
{\bf Proof}: We have $d\xi=Im\varphi'$. On each $S_i$ the action is non-regular, thus $\varphi'=0$ on $S_i$. In particular $d\xi|_{S_i}=0$, i.e $\xi$ is a constant $\xi_i$ on $S_i$.

Let $0 \neq e_i$ be an element in the Lie algebra of $T_i$. Then the flow v.field $X_i$ of $e_i$ vanishes along $S_i$. Thus $d\mu(e_i)=i_{X_i}\omega=0$ along $S_i$. So in particular $\mu(e_i)=\mu_i=const$. So the image of $\mu$ on $S_i$ lives on a hyperplane in the dual Lie algebra ${\cal G}^{\ast}$ of $T$.

From all this we deduce that the image of $(\mu,\xi)$ on $S$ lives on a finite union $H$ of $(n-2)$-planes in $\mathbb{R}^n$.       Q.E.D. 
\section{Group actions and coassociative submanifolds}
Let $M^7$ be a 7-manifold with a $G_2$-form $\varphi$. This means that for each point $m \in M$ there is an isomorphism $\sigma$ between the tangent space $T_mM$ and $\mathbb{R}^7$ s.t. $\sigma^{\ast}(\varphi_0)= \varphi$, there $\varphi_0$ is the standard $G_2$ form on $\mathbb{R}^7$ (see \cite{LM}). Since the group $G_2$ preserves the Cayley product on $\mathbb{R}^8= \mathbb{R}^7 \oplus \mathbb{R}$, then the bundle $TM \oplus \mathbb{R}$ over $M$ acquires a structure of an algebra, isomorphic to Cayley numbers (see \cite{LM}).

We will assume that $\varphi$ is closed and co-closed, hence parallel and the Holonomy of $M$ is contained in the group $G_2$. The 4-form $\ast \varphi $ is a calibration and a calibrated submanifold $L$ is called a coassociative submanifold. This is equivalent to $\varphi|_L=0$. We will use group actions to construct coassociative submanifolds on $M$. We assume that 
$b_1(M)=0$. We will treat 3 cases :

1) A 3-torus action. Let $v_i$ be a basis of the Lie Algebra of $T^3$ and $X_i$ are
the corresponding (commuting) flow vector fields on $M$. Then we define the 1-forms 
$\sigma_1= i_{X_2}i_{X_3} \varphi$, $\sigma_2 = i_{X_3}i_{X_1}\varphi$,
$\sigma_3= 
i_{X_1}i_{X_2}\varphi$. As in section 2, one can easily show that $\sigma_i$ are closed, $T^3$-invariant 1-forms. Hence $\sigma_i=df_i$ for some $T^3$-invariant functions $f_i$. Since $f_1$ is $T^3$-invariant we have \[0= df_1(X_1)=
\sigma_1(X_1)= \varphi(X_1,X_2,X_3)\]
Consider now a level set $L=(f_i=const)$. Suppose some point $m \in L$ is a regular point of $T^3$-action. Since $\varphi(X_1,X_2,X_3)=0$ one easily sees that $\sigma_i$ are linearly independent at $m$ and hence $m$ is a smooth point of $L$. Also $L$ is $T^3$-invariant, hence $X_i(m)$ are in the tangent bundle of $L$ at $m$. Hence one easily deduces that $\varphi |_L=0$, i.e. $L$ is coassociative.

2) A 2-torus action. We have vector fields $X_1$ and $X_2$ and a 1-form $\sigma= i_{X_1}i_{X_2}\varphi$. Once again $\sigma=df$ for a ($T^2$-invariant) f. Consider a level set $N=(f=const)$. Suppose $T^2$ acts freely on $N$ and consider the quotient $M_{red}=N/T^2$, which we call a $G_2$-reduction. We have a projection $\pi: N \mapsto M_{red}$. Consider the (closed) 2-forms $\eta_i=i_{X_i}\varphi $ and $\omega= i_{X_1}i_{X_2} (\ast \varphi)$. One can easily show that there are unique, closed 2-forms $\eta_i'$ and $\omega'$ on $M_{red}$ s.t. $\pi^{\ast}(\eta_i')=\eta_i$ and $\pi^{\ast}(\omega')=\omega$. Let $L'$ be a 2-submanifold of $M_{red}$ and $L= \pi^{-1}(L')$. Then obviously $L$ is coassociative iff $\eta_i'|_{L'}=0$. We will reformulate this condition as a pseudoholomorphic condition on $L'$.

Consider a bundle $V= span(X_1,X_2)$ over $N$. Pick $n \in N$ and let $e_1,e_2$ be an orthonormal basis of $V$ at $n$, compatible with the orientation, given by $X_1$ and $X_2$. Let $e_3= e_1 \times e_2$ (here $\times$ is the Cayley product). $e_3$ doesn't depend on the choice of $e_1$ and $e_2$ and thus gives rise to a section of $TM$ over $N$. Consider the bundle $W$ over $N$, which is the orthogonal complement of $V \oplus (e_3)$ in $TM$. Then one easily verifies that the tangent bundle of $N$ is $W \oplus V$ and the quotient of $W$ by $T$-action can be viewed as a tangent bundle to $M_{red}$. Let again $n \in N$ and let $J_i$ be a right Cayley multiplication by $e_i$. Then $J_i$ preserve the fiber $W_n$ of $W$ at $n$ and they give complex structures on $W_n$, which form a HyperKahler package. Also $J_3$ gives rise an almost complex structure on $M_{red}$. Let $\omega_i$ be the corresponding symplectic forms on $W_n$. Then one easily verifies that $\omega$ is proportional to $\omega_3$ on $W_n$ and $span(\eta_1|_{W_n},\eta_2|_{W_n})=span(\omega_1,\omega_2)$ in the space of 2-forms on $W_n$. From all this linear algebra we get that $\omega'$ is a symplectic form on $M_{red}$ and $J_3$ is a compatible almost complex structure. Also for a 2-submanifold $L'$ of $M_{red}$ the conditions $\eta_i'|_{L'}=0$ are equivalent to $L'$ being $J_3$-holomorphic.

3) An $SO(3)$ action. We will also assume that the action is not regular in at least one point.
Let $e_1,e_2,e_3$ be the o.n. basis of the Lie Algebra of $SO(3)$. We have the following relations : \[ [e_1,e_2]=e_3 ~ , ~ [e_2,e_3]=e_1 ~ , ~ [e_3,e_1]= e_2 \]   
Let $X_i$ be the corresponding vector fields on $M$. Let $\sigma= i_{X_1}i_{X_2}\varphi$. Then \[d\sigma= {\cal L}_{X_1}(i_{X_2}\varphi)= i_{[X_1,X_2]}\varphi=i_{X_3}\varphi \]
Let $f= \varphi(X_1,X_2,X_3)$. Then \[df= {\cal L}_{X_3}\sigma - 
i_{X_3}d\sigma\] We easily deduce that both terms are $0$, hence $df=0$. Also $f=0$ in at least 1 point. Hence $f=0$. So we might hope to find coassociative submanifolds, invariant under $SO(3)$-action.  

Let $\alpha= i_{X_1}i_{X_2}i_{X_3} (\ast \varphi)$. Then using arguments as before we can show that $\alpha$ is a closed, $SO(3)$-invariant 1-form. So $\alpha=dg$ for an $SO(3)$-invariant function $g$. Let $v= \nabla g $. Let $m$ be a regular point of the action. Then $v \neq 0$ at $m$. Also the scalar product of $v$ and $X_i$ is $0$, so $v$ and $X_i$ are linearly independent and span a 4-dimensional space $W$. Using some Cayley algebra one can easily show that $W$ is a coassociative subspace of $TM$. 

We assume that $SO(3)$ acts freely on the space $M'$ of regular points of the action. Also the complement $M-M'$ corresponds precisely to the critical points of $g$. Let $l$ be a non-constant trajectory of the gradient flow of $g$. Then $l$ is contained in $M'$.
Let $L =SO(3) \times l$. Then $L$ is coassociative. Also trajectories of the gradient flow are embedded 1-submanifolds, $g$ is $SO(3)$-invariant and increases on the trajectories. From all this we deduce that $L$ is an embedded submanifold. Thus $M'$ is covered by a family of non-intersecting coassociative submanifolds, diffeomorphic to $SO(3) \times \mathbb{R}$.

We can't in general say anything about the set of non-regular points. We will do this in one example in section 4. 
\section{Examples}
In this section we will give a number of examples, there results of the two previous sections are applicable.
\subsection{$\mathbb{C}^n$}
There is a $T^{n-1}$ action on $\mathbb{C}^n$ given by \[(e^{i \theta_1},\ldots,e^{i \theta_{n-1}})\cdot(z_1,\ldots,z_n)= (e^{i\theta_1}z_1,\ldots,
 e^{i\theta_{n-1}}z_{n-1}, e^{-i(\theta_1 + \ldots + \theta_{n-1})}z_n) \] 
The moment map for this action is $\mu=(\mu_1,\ldots,\mu_{n-1})$ with 
$\mu_i= |z_i|^2-|z_n|^2$. The 1-form $\varphi'$, defined in the proof of Theorem 1, is
$\varphi'= i^{n-1}\Sigma dz_i ( z_1\cdots \stackrel{\wedge}{z_i}\cdots z_n)= d (i^{n-1}z_1 \cdots z_n)$. So the fibration is given by
\[ |z_i|^2-|z_n|^2=c_i ~ , ~ Im(i^{n-1} z_1 \cdots z_n)=c_n \]
This is a classical example of Harvey and Lawson (see \cite{HL}).
\subsection{$K(N)$ and the Calabi construction}
Consider $\mathbb{C}^n$ and a $\mathbb{Z}_n$-action on it with $k \in
 \mathbb{Z}_n $ acts by multiplication by $e^{2\pi ki/n}$. Then the quotient has a resolution of singularities, which is a total space of $\gamma ^{\otimes n}$, there $\gamma$ is the universal line bundle over $\mathbb{C}P^{n-1}$, i.e. the resolution of singularities is the total space of the canonical bundle over $\mathbb{C}P^{n-1}$. 

Let $K(N)$ be a total space of a canonical bundle of a complex manifold $N$ and $\pi : K(N) \mapsto N$ be a projection. There is a canonical $(n,0)$-form $\rho$ on $K(N)$ defined by $\rho(a)(v_1,\ldots,v_n)=a(\pi_{\ast}(v_1),\ldots,\pi_{\ast}(v_n))$, $a \in K(N)$.
The form $\varphi = d\rho$ is a holomorphic volume form on $K(N)$. If $z_1, \ldots,z_n$ are local coordinates on $N$ then
 $(z_1,\ldots,z_n,y=dz_1\wedge \ldots \wedge dz_n)$ are coordinates on $K(N)$ and $\varphi= dz_1\wedge \ldots \wedge dz_n\wedge dy$.

$\mathbb{C}P^{n-1}$ is a Kahler-Einstein manifold with positive scalar curvature. If $N$ is a K-E manifold with positive scalar curvature then $K(N)$ has a Ricci-flat Kahler metric on it (see \cite{Sol}, p.108). The metric is constructed as follows : The connection on $K(N)$ induces a horizontal distribution for the projection $\pi$, with a corresponding splitting of the tangent bundle of $K(N)$ into horizontal and vertical distributions. We can view the horizontal space at each point $m \in K(N)$ as a tangent space to $N$ at $\pi(m)$. Let $r^2 : K(N) \mapsto 
\mathbb{R}_+$ be the square of the length of an element in $K(N)$ and $u: \mathbb{R}_+ \mapsto \mathbb{R}_+$ be a positive function with a positive first derivative. We define the metric $\omega_u$ on $K(N)$ as follows: We put the horizontal and the vertical distributions to be orthogonal. On the horizontal distribution we define the metric to be $u(r^2)\pi^{\ast}(\omega)$ and on the vertical distribution we define it to be $t^{-1}u'(r^2)\omega'$. Here $\omega$ is the Kahler-Einstein metric on $N$, $t$ is it's scalar curvature and $\omega'$ is the induced metric on linear fibers of $\pi $. The Kahler-Einstein condition ensures that the corresponding 2-form $\omega_u$ defining this metric on $K(N)$ is closed, i.e. the metric is Kahler. If we take $u(r^2)=(tr^2+l)^{\frac{1}{n+1}}$ for some constant $l$  (see \cite{Sol}, p.109), then $\omega_u$ is Ricci-flat.

For $\mathbb{C}P^{n-1}$ the Ricci-flat metric on $K(\mathbb{C}P^n)$ has a Kahler potential $f$ outside of the zero section and $f$ is a function $f=h(r^2)$, there $r^2= \Sigma|z_i|^2$ on $(\mathbb{C}^n- 0)/\mathbb{Z}_n$. For instance then $n=2$ we have the Eguchi-Hanson potential $h(x)= \sqrt{x^2+1} + logx- log(\sqrt{x^2+1}+1)$. Also the metric is asymptotic to the Euclidean metric on $\mathbb{C}^n/\mathbb{Z}_n$ at infinity.

We have an $(n-1)$-torus action on $\mathbb{C}^n$ as in the first example, and this action commutes with the $\mathbb{Z}_n$-action, hence it induces an action on $K(\mathbb{C}P^{n-1})$. This action preserves the Calabi-Yau structure on $K(\mathbb{C}P^{n-1})$, and hence Theorem 2 applies. To compute the moment map, we note that $\omega_u= i \overline{\partial} \partial f= d(i \partial f)$. Now the $T^{n-1}$-action preserves $f$ and so it preserves $\partial f$. Let $v \in {\cal G}$ and $X_v$ be the vector field on $K(\mathbb{C}P^{n-1})$ as before. Then
\[0 = {\cal L}_{X_v}(i \partial f)= i_{X_v}\omega_u + d(i \partial f (X_v)) \]
Now $ i \partial f (X_v) = i (df(X_V)- i df(JX_v))$. Also $df(X_v)=0$, so 
$i \partial f (X_v)= df(JX_v)=h'(r^2)\cdot dr^2(JX_v)$. If $v=\partial_{\theta_i}$ then $dr^2(Jx_v)= |z_n|^2-|z_i|^2$.

So the moment map $\mu=(\mu_1, \ldots,\mu_{n-1})$ 
satisfies $d\mu_i= i_{X_i}\omega_u=d(h'(r^2) \cdot (|z_i|^2 -|z_n|^2))$.
So $\mu_i=h'(r^2)(|z_i|^2-|z_n|)^2$.
 By similar reasoning the function $\xi$ is given, as in the previous example, by $\xi= Im(i^{n-1} z_1 \cdots z_n)$. So SLag fibration on $K(\mathbb{C}P^n)$ is given by 
\begin{equation} 
h'(r^2) \cdot (|z_i|^2-|z_n|^2) = c_i ~ , ~ Im(i^{n-1} z_1 \cdots z_n) = c_n 
\end{equation}
Also $h'(r^2)$ converges to 1 at infinity, so this fibration is asymptotic at infinity to the fibration on $\mathbb{C}^n$ in the previous example.

We will now make the following general observation : Let $L$ be an oriented Lagrangian submanifold of a K-E manifold $N$. We endow $K(N)$ with a Kahler metric $\omega_u$ as above for any choice of a function $u$. For any point $l \in L$ there is
a unique element $\kappa_l$ in the fiber of $K(N)$ over $l$ which restricts to the volume form on $L$. Various $\kappa_l$ give rise to a section $\kappa$ of $K(N)$ over $L$.
 For any $\lambda \in \mathbb{R}$ we define a submanifold $L_{\lambda} \subset K(N)$ by \[ L_{\lambda}=((l,\mu)|l \in L ~ , ~ \mu=(a+i\lambda)\kappa_l ~ , ~a \in \mathbb{R}) \] 
We have the following:
\begin{lem}
$L$ is a minimal Lagrangian submanifold of $N$ iff any of $L_{\lambda}$ is a Special Lagrangian submanifold of $N(K)$
\end{lem} 
{\bf Proof :} First we note that $L_{\lambda}$ are Special, i.e.
$Im\varphi|_{L_{\lambda}}=0$. Indeed one easily verifies that $Im\rho|_{L_{\lambda}}=\lambda \pi^{\ast}(\kappa)$ and hence $Im\varphi|_{L_{\lambda}}=\lambda
\pi^{\ast}(d\kappa)=0$. 

We now prove that $L_{\lambda}$ is Lagrangian iff $L$ is minimal. Let $m$ be a point on $L_{\lambda}$, $l = \pi(m)$ and $m = (a+i\lambda)\kappa_l$. The tangent space of $L_{\lambda}$ at $m$ is spanned by $\kappa_l$ (viewed as a vertical vector in $T_mK(N)$) and vectors $(e+(a+i\lambda)
\nabla_e\kappa)$. Here $e$ is any tangent vector to $L$ at $l$ (viewed as an element of the horizontal
distribution of $T_mK(N)$) and $(a+i\lambda)\nabla_e \kappa$ lives in the vertical distribution of $T_mK(N)$.
To compute $\nabla_e \kappa$ take an orthonormal frame $(v_j)$ of $T_lL$ and extend it to an orthonormal frame of $L$ in a neighbourhood $U$ of $l$ in $L$ s.t. $\nabla^L v_i=0$ at $l$ (here $\nabla^L$ is the Levi-Civita connection of $L$). We get that \[ \nabla_e \kappa = \kappa \cdot \nabla_e \kappa(v_1,\ldots,v_n)= \kappa (e(\kappa(v_1,\ldots,v_n))- \Sigma\kappa(v_1,\ldots,\nabla_e v_j,\ldots,v_n)) \]
Now $e(\kappa(v_1,\ldots,v_n))=0$. Also clearly \[\kappa(v_1,\ldots,\nabla_e v_j,\ldots,v_n)=i<\nabla_e v_j,Jv_j>=i<\nabla_{v_j}e,Jv_j>=-i<e,J(\nabla_{v_j}v_j)> \]   
Here $J$ is the complex structure on $N$. Thus we get that 
\[(a + i\lambda) \nabla_e \kappa = i(a+i \lambda)(Jh \cdot  e)\kappa_l \]
Here $h= \Sigma \nabla_{v_j}v_j$ is the trace of the second fundamental form
of $L$. From this one easily deduces that $L_{\lambda}$ is Lagrangian iff
$h=0$, i.e. $L$ is minimal.    Q.E.D.

We will now investigate toric K-E manifolds. For recent results  and examples we refer the reader to \cite{Bat} and \cite{Tian}. We begin with the following lemma.
\begin{lem}
Let $(M^{2n},\omega)$ be a compact symplectic manifold and $g$ some Riemannian metric on $M$. Suppose that we have an effective Hamiltonian n-torus action on $M$, which preserves $g$. Then there is a regular orbit of the action, which is a minimal submanifold with respect to $g$.
\end{lem}
Here by regular orbits we mean orbits with a finite stabilizer.

{\bf Proof:} We have a moment map $\mu$ and smooth orbits are levels set of the moment map. For an orbit $L$ to be a minimal submanifold, it is obviously necessary to be a critical point of the volume functional on the orbits. We note that it is also sufficient. Indeed let $v$ be any element of the Lie Algebra $\mathcal{G}$ of the torus $T^n$. Then $\mu(v)$ is $T^n$-invariant, and so is the gradient
$\nabla \mu(v)$. Also this gradient is orthogonal to the orbits. Consider now this gradient flow. It commutes with $T^n$-action, hence it sends orbits to orbits. Since $L$ is critical for the volume functional, we get from the first variation formula $\int_{L}h \cdot \nabla \mu(v)=0$. Here $h$ is a trace of the second fundamental form of $L$. But both $h$ and $\nabla \mu(v)$ are $T^n$-invariant, hence we are integrating a constant. So $h \cdot \nabla \mu(v)=0$ pointwise. Now $v$ was arbitrary, hence $h=0$.  

Finally we note that at least one orbit is critical for the volume functional on the orbits. We use the following easy
\begin{lem}
Let $L$ be a orbit with a positive dimensional stabilizer $T' \subset T$ and $x \in L$. Then for any $\epsilon > 0$ there is a neighbourhood $U$ of $x$ s.t. any orbit passing through $U$ has volume $< \epsilon$.
\end{lem}
Indeed we can take a (unit) vector $e_1$ in the Lie Algebra of $T'$. Then the corresponding flow vector field $X_1$ vanishes along $L$. Extend $e_1$ to an o.n. basis $e_2,\ldots,e_n$ of the Lie algebra of $T$. The vector fields $X_i$ will have uniformly bounded lengths. We choose a neighbourhood $U$ of $x$ in which $X_1$ has sufficiently small length and it is clear that volumes of orbits through $U$ will be sufficiently small.

So now we try maximize volume among regular orbits. Let $L_i$ be a sequence of orbits, whose volume goes to supremum $s$ of volumes of regular orbits. Then by the previous lemma it is clear that a limiting orbit (of some subsequence) $L$ is regular. Now the differential of the moment map on $M$ is surjective along $L$ and $L_i$ are level sets of $\mu$. It is clear that $L_i \mapsto L$ as manifolds and hence the volume of $L$ is $s$ and we are done.            Q.E.D.
 
On $\mathbb{C}P^{n-1}$ we have the following $T^{n-1}$-invariant minimal Lagrangian torus \[L=((z_1, \ldots,z_n)| |z_i|=|z_j|)\] and corresponding SLag submanifolds $L_{\lambda}$ in $K(\mathbb{C}P^{n-1})$. Those submanifolds are invariant under
our $T^{n-1}$-action and they are in the moduli-space $(2)$ we constructed (in fact $L_0$ is a submanifold, which corresponds to $c_i=0$ in equation (2)). Moreover any other element in our moduli-space is asymptotic at infinity to $L_0$. By that we mean the following : Let $B$ be the unit ball of $K(\mathbb{C}P^{n-1})$ with respect to $r^2$. $L_0$ is invariant under scaling of $K(\mathbb{C}P^{n-1})$ by a real number. If $L$ is another element in our moduli-space then scaling of $L$ by $k \in \mathbb{R}$ as $k$ goes to infinity converges in $C^{\infty}$ to $L_0$ on compact subsets of $B - \mathbb{C}P^{n-1}$. It turns out that analogous situation holds for any toric K-E manifold N with positive scalar curvature.

Suppose we have an effective, structure-preserving $T^n$-action on $N$.
We make the following definition: We have a $T^n$ action on $N$ and this
action induces a $T^n$-action on
$K(N)$. Let $\mathcal{G}$ be the Lie algebra of $T$, $ v \in \mathcal{G}$,
$X_v$ be
the flow vector field on $N$ and $X_v'$ be the flow vector field on $K(N)$.
So $\pi_{\ast}(X_v')=X_v$. Let $l \in N$ and $m \in K_l=\pi^{-1}(l)$.
Let $R(m)$ be the vertical part of $X_v'$ at $m$. Since $R(m)$ is vertical, it
can be viewed as an element of $K_l$. The correspondence $m \mapsto R(m)$ is a
linear correspondence on $K_l$. Hence there is a complex number $\sigma_l(v)$
s.t. $R(m)=\sigma_l(v)m$. At a regular point $l$ of the $T$-action $\sigma_l(v)$ can also be found in a following way : Take any unit length element $\xi \in K_l$. Extend $\xi$ along the orbit of $X_v$ to be invariant under the flow of $X_v$. Then one easily computes that $\sigma_l(v)= \nabla_{X_v}\xi \cdot \xi$. Since the flow of $X_v$ is given by holomorphic isometries, $\xi$ has unit length. Hence  $\sigma_l(v)$ is purely imaginary. Also $\sigma_l(v)$ is linear in $v$ (because $R(m)$ is given by the vertical part of the differential of the $T$-action at $m$, and this differential is a linear map from $\mathcal{G}$ to $T_mK(N)$). Hence $i \sigma$ can be viewed as a map from $N$ to the dual Lie algebra $\mathcal{G}^{\ast}$. This map is $T$-invariant.
 
\begin{lem}
For a regular orbit $L$ of the $T^n$-action, $L$ is minimal iff $\sigma|_L = 0$
\end{lem}
{\bf Proof:} The section $\kappa$ of $K(N)$ over $L$ we defined in Lemma 1 is $T^n$-invariant. Also we computed that $\nabla_{X_v}\kappa \cdot \kappa = i (Jh \cdot X_v)$ for any $v \in {\cal G}$. From this the lemma follows     Q.E.D.

Let $t$ be the scalar curvature of $N$.
\begin{lem}
The map $\mu = -it^{-1} \sigma$ is a moment map for the action. The zero set of $\mu$ is precisely 1 regular orbit $L$.
\end{lem}
{\bf Proof} Let $v \in \mathcal{G}$. We need to show that $d(-it^{-1}\sigma(v))= i_{X_v}\omega$. We will do it at a smooth point $p$ of the action. Choose any unit length element $\xi$ of $K(N)$ over $p$. We can extend $\xi$ to be a local unit length section, invariant under the $X_v$-flow. We have a connection 1-form $\eta(u)= \nabla_u \xi \cdot \xi$. Then $\eta$ is invariant under the $X_v$-flow and the K-E condition says that $i d\eta = t \omega$. 
So \[ 0 = \mathcal{L}_{X_v}\eta= d(i_{X_v}\eta) + i_{X_v}d\eta= d\sigma(v) - it (i_{X_v}\omega) \]   
So $\mu$ is a moment map. By Lemma 2 one of the regular orbits $L$ is minimal, hence it lies on the zero set of $\mu$ by Lemma 4. Obviously this orbit is isolated in the zero set of $\mu$. Now by Atiyah's result (see \cite{At}), the zero set of the moment map is connected, hence it must coincide with $L$ and we are done.
   Q.E.D. 
\begin{lem} 
The map $\mu'= u\pi^{-1}(\mu)$ is a moment map for the $T$-action on
$K(N)$.
\end{lem}
{\bf Proof:} Let $v \in \cal{T}$. We need to prove that $d\mu'(v)=i_{X_v'}\omega_u$.
 
We will study $\omega_u$ in more detail (see \cite{Sol}). Let $m \in N$ be a regular point for the $T^n$-action and $\xi$ a unit length element of $K(N)$ over $m$. We can extend $\xi$ to be a local unit length section of $K(N)$, invariant under the flow of $X_v$. $\xi$ gives rise to a connection 1-form $\psi$ for the connection on $K(N)$  and the Einstein condition tells that $id\psi= t \omega$. The section $\xi$ defines a complex coordinate $a$ on $K(N)$, which is invariant under the $X_v'$-flow. Also the form $b=da + a\pi^{\ast}\psi$ vanishes on the horizontal distribution (see \cite{Sol}, p. 108). We have $r^2=a\overline{a}$ and $u= u(r^2)$. Also the Kahler form $\omega_u$ on $K(N)$ is given by \[\omega_u=u\pi^{\ast}\omega -it^{-1}u'b\wedge \overline{b} \]
One directly verifies that $\omega_u = d\eta$ for $\eta= it^{-1}u \pi^{\ast} \psi - it^{-1}\frac{ud\overline{a}}{\overline{a}}$. By our construction the flow of $X_v'$ leaves $\eta$ invariant. So \[0 = {\cal L}_{X_v'}\eta= i_{X_v'}d\eta + d(i_{X_v'}\eta)= i_{X_v'}\omega_u + d(it^{-1}u\psi(X_v))= i_{X_v'}\omega- d(\mu'(v))\] Here we used the fact that $d\overline{a}(X_v')= 0$ and $\psi(X_v)= \sigma(v)$. So $\mu'$ is a moment map and we are done.   Q.E.D.

The torus action on $N$ induces an action on $K(N)$ and by Theorem 2 we have a SLag fibration on $K(N)$. We want to investigate the asymptotic behavior of the fibers at infinity. We will assume that the function $u=u(r^2)$, used to define the metric $\omega_u$ on $K(N)$ goes to infinity as $r^2$ goes to infinity (this holds e.g. for $u$ defining the Ricci-flat metric).

\begin{thm} 
$L_0 \subset K(N)$ is a fiber of the fibration arising from
Theorem 2. Moreover, any other fiber is asymptotic to it at infinity.
\end{thm} 
{\bf Proof:} Let $e_1,\ldots,e_n$ be a basis for $\mathcal{G}$. Let $X_i$ be the flow fields on $N$ and $X'_i$ be the flow fields on $K(N)$. Then $\pi_{\ast}(X'_i)=X_i$. 
Let $\rho$ be an $(n,0)$ form on $K(N)$ as before and $\varphi=d\rho$. Then $\rho$ is $T^n$-invariant and we can prove, as in Theorem 1, that $i_{X'_1}\ldots i_{X'_n}\varphi=d(\rho(X'_1,\ldots,X'_n))$. Also for $\xi \in K(N)$, $\rho(\xi)(X'_1,\ldots,X'_n)=\xi(X_1(\pi(\xi)),\ldots,X_n(\pi(\xi)))$.
We also have a moment map $\mu'=
-it^{-1}u\pi^{-1}(\sigma)=u\pi^{-1}(\mu)$. Thus our SLag fibers will be
given by equations
\[u(|\xi|^2) \mu_i(\pi(\xi))=c_i~ , ~ Im\xi(X_1(\pi(\xi)),\ldots,X_n(\pi(\xi)))=c_n \]
Here $\xi \in K(N)$. The fiber $L_0$ corresponds to $c_j=0$.
The asymptotic behavior of the fibers follows immediately from this
formula. Q.E.D.

We can ask a similar question in general, thus giving a non-compact analog of the SYZ conjecture (see \cite{SYZ}) : Let $N$ be a K-E manifold with positive scalar curvature. When is it true that $N$ has a minimal Lagrangian submanifold $L$ and $K(N)$ is fibered by SLag subvarieties with fibers asymptotic to $L_0$ at infinity ?   
\subsection{Resolutions of singularities and (Quasi)ALE spaces}
Suppose that a finite subgroup $G$ of $SU(n)$ acts on $\mathbb{C}^n$ and we have a crepant resolution of singularities $M$. D. Joyce has recently constructed a (Quasi)ALE Ricci-flat Kahler metric $\omega$ on $M$ (see \cite{Joy1} and \cite{Joy2}). 

Suppose that $G$ is Abelian. Then there is an orthonormal basis of $\mathbb{C}^n$, in which the action of $G$ is given by diagonal matrices. 
Consider now the $T^{n-1}$-action on $\mathbb{C}^n$ as in the first example (4.1). 
 This action commutes with the $G$-action, hence it induces an action on $M$. Also by the uniqueness property of Joyce's construction, this action preserves $\omega$. Hence Theorem 2 applies, and we have a Special Lagrangian fibration on $M$.
\subsection{Coassociative submanifolds} 
Robert Bryant and Simon Salamon have constructed in \cite{BS} some examples of complete $G_2$ metrics. Some examples are on total space of a spin bundle over a 3-dimensional space form. Others are on total space $\Lambda^2_-$ of anti-self-dual 2-forms over a self-dual Einstein 4-manifold. Those 3 and 4-manifolds admit isometric actions by 2-tori and by $SO(3)$, and those actions induce structure-preserving actions on the corresponding $G_2$-manifolds. We will treat 1 example in detail- the total space of a spin bundle over $S^3$. 

The spin bundle is a bundle $V=TS^3 \oplus \mathbb{R}-$ the direct sum of the tangent bundle of $S^3$ with a trivial bundle. Now $S^3$ can be viewed as a unit sphere of quaternions. There is an $S^3$ action on itself, given by $q(p)=qpq^{-1}$. Here the multiplication is a quaternionic multiplication. Obviously this action becomes an $S^3/ \pm1 = SO(3)$-action. 

In \cite{BS} a $G_2$-structure was constructed on the total space of $V$. We won't reproduce the details of the construction but only mention that the fibers of the projection of $V$ over $S^3$ are coassociative. We will look for coassociative submanifolds, invariant under $SO(3)$ action.

The points $\pm 1$ are fixed by $SO(3)$ action and the fibers over these points are $SO(3)$-invariant coassociative submanifolds $L_{\pm 1}$. Take now any point $m \in (S^3 - \pm 1)$. Then the stabilizer of $SO(3)$ action on $m$ is a circle. Let $N_m$ be the orthogonal complement in the tangent space $T_m S^3$ to the orbit of $SO(3)$-action. Then $W=N_m \oplus \mathbb{R}$ is a sub-bundle of $V$ over $(S^3 - \pm 1)$. $W$ is invariant under $SO(3)$-action. Let $A_m$ be the orbit of $m$ under $SO(3)$-action ($A_m$ is diffeomorphic to $S^2$) and $L_m$ be the total space of $W$ over $A_m$.
 Then one can easily show that $L_m$ is a coassociative submanifold, invariant under $SO(3)$-action. Also the union of all $L_m$ and of $L_{\pm1}$ is precisely the set of non-regular points of the action. Also $SO(3)$ acts freely on the set of regular points of the action.
By the results of section 3, the set of regular points is covered by a family of non-intersecting coassociative submanifolds. So the whole $V$ is covered by non-intersecting, coassociative submanifolds. Those are a 3-dimensional family of submanifolds, diffeomorphic to $SO(3) \times \mathbb{R}$, a 1-dimensional family of submanifolds, diffeomorphic to $S^2 \times \mathbb{R}^2$ and 2 submanifolds, diffeomorphic to $\mathbb{R}^4$.

Massachusetts Institute of Technology

E-Mail : egold@math.mit.edu

\end{document}